\documentclass[AMA,STIX1COL]{WileyNJD-v2}

\usepackage{amsmath}
\usepackage{bm}
\usepackage{tikz}
\usetikzlibrary{positioning}
\usepackage{standalone}

\articletype{Research Article}%

\begin{document}

\title{A note on the conservation properties of the generalized-$\alpha$ method}

\author{DeAnna S. Gilchrist}

\author{John A. Evans}

\authormark{GILCHRIST \textsc{et al}}

\address{\orgdiv{Ann and H.J. Smead Department of Aerospace Engineering Sciences}, \orgname{The University of Colorado Boulder}, \orgaddress{\city{Boulder}, \state{Colorado}, \country{USA}}}

\corres{John A. Evans, Ann and H.J. Smead Department of Aerospace Engineering Sciences, The University of Colorado Boulder, Boulder, CO, USA.\\ \email{john.a.evans@colorado.edu}}

\abstract[Abstract]{We show that the second-order accurate generalized-$\alpha$ method on a uniform temporal mesh may be viewed as an implicit midpoint method on a shifted temporal mesh.  With this insight, we demonstrate generalized-$\alpha$ time integration of a finite element spatial discretization of a conservation law system results in a fully-discrete method admitting discrete balance laws when (i) the time integration is second-order accurate, (ii) a uniform temporal mesh is employed, (iii) the spatial discretization is conservative, and (iv) conservation variables are discretized.}


\maketitle
\section{Introduction}\label{sec1}

The generalized-$\alpha$ method is a family of time integration schemes that, for a particular choice of method parameters, is second-order accurate, is unconditionally stable, and exhibits an optimal combination of accuracy in the low-frequency range and damping in the high-frequency range.  The generalized-$\alpha$ method was first introduced for second-order initial value problems by Chung and Hulbert \cite{chung1994family}, and it was later extended to first-order initial value problems by Jansen, Whiting, and Hulbert \cite{jansen2000generalized}.  The generalized-$\alpha$ method is typically combined with a finite element spatial discretization in order to arrive at a fully-discrete method for the numerical solution of partial differential equations.  This is a particularly popular approach for structural mechanics applications \cite{hulbert1996explicit,kuhl1999energy,kuhl1999generalized,arnold2007convergence,erlicher2002analysis}, though it is often used for fluid mechanics \cite{bazilevs2007variational,gomez2010isogeometric,modirkhazeni2016algebraic,bayram2020variational,codoni2021stabilized,liu2021note}, fluid-structure interaction \cite{dettmer2006computational,bazilevs2008isogeometric}, and magnetohydrodynamics \cite{gleason2022divergence} applications as well.  In the context of fluid mechanics, the generalized-$\alpha$ method is used to time integrate finite element spatial discretizations of mass, momentum, and energy differential conservation laws.  However, this results in a fully-discrete method that does not admit discrete balance laws for mass, momentum, and energy with respect to the temporal mesh, even if the underlying spatial discretization is conservative.  Fortunately, we show in this note that the resulting fully-discrete method does admit discrete balance laws for mass, momentum, and energy with respect to a shifted temporal mesh if the temporal mesh is uniform and if the parameters of the generalized-$\alpha$ method are chosen so that it is second-order accurate.  To arrive at this result, we invoke a new interpretation of the second-order accurate generalized-$\alpha$ method.  Namely, it can be interpreted as an implicit midpoint method on a shifted temporal mesh when the temporal mesh is uniform.

An outline of this short communication is as follows.  In Section \ref{sec:gen_alpha}, we show how the generalized-$\alpha$ method for first-order initial-value problems can be viewed as an implicit midpoint method on a shifted temporal mesh when it is second-order accurate and the temporal mesh is uniform.  In Section \ref{sec:advection_diffusion}, we use this knowledge to show that application of second-order accurate generalized-$\alpha$ time integration to a conservative stabilized or unstabilized Galerkin discretization of a model advection-diffusion problem yields a fully-discrete method harboring a discrete balance law when the temporal mesh is uniform.  In Section \ref{sec:conservation_laws}, we show the same is true for general systems of differential conservation laws provided the conservation variables are themselves discretized, and in Section \ref{sec:nonconservative variables}, we show how to modify the generalized-$\alpha$ method to arrive at a conservative fully-discrete method when nonconservation variables are discretized instead.  Finally, in Secion \ref{sec:conclusion}, we provide concluding remarks.

\section{An alternative form of the generalized-$\alpha$ method}\label{sec:gen_alpha}

Consider the following first-order initial-value problem: Find $\mathrm{U}: [0,\infty) \rightarrow \mathbb{R}^m$ such that
\begin{equation}
\mathrm{R}\left(\dot{\mathrm{U}}(t),\mathrm{U}(t),t\right) = \mathrm{0}
\end{equation}
for all $t \in (0,\infty)$ and
\begin{equation}
\mathrm{U}(0) = \mathrm{U}_0
\end{equation}
where $m \in \mathbb{N}$ is the size of the solution vector $\mathrm{U}$, $\dot{\mathrm{U}}$ is the time derivative of $\mathrm{U}$, $\mathrm{U}_0 \in \mathbb{R}^m$ is the initial condition of $\mathrm{U}$, and $\mathrm{R}: \mathbb{R}^m \times \mathbb{R}^m \times (0,\infty) \rightarrow \mathbb{R}^m$ encodes the ordinary differential equations associated with the initial-value problem.  In the generalized-$\alpha$ method, the solution vector $\mathrm{U}$ is approximated on a temporal mesh $\left\{ t_n \right\}_{n \in \mathbb{N}}$ of increasing times with $t_1 = 0$.  In particular, given the approximations $\dot{\mathrm{U}}_n$ and $\mathrm{U}_n$ of $\dot{\mathrm{U}}$ and $\mathrm{U}$ at the $n^{\text{th}}$ time $t_n$, the generalized-$\alpha$ method involves solving the following algebraic system of equations for the approximations $\dot{\mathrm{U}}_{n+1}$ and $\mathrm{U}_{n+1}$ of $\dot{\mathrm{U}}$ and $\mathrm{U}$ at the $(n+1)^{\text{st}}$ time $t_{n+1}$\cite{jansen2000generalized}:
\begin{align}
\mathrm{R}\left(\dot{\mathrm{U}}_{n+\alpha_m},\mathrm{U}_{n+\alpha_f},t_{n+\alpha_f}\right) = \mathrm{0} \label{eq:gen_alpha}
\end{align}
and
\begin{align}
\mathrm{U}_{n+1} = \mathrm{U}_{n} + \Delta t_n \left( (1-\gamma) \dot{\mathrm{U}}_{n} + \gamma \dot{\mathrm{U}}_{n+1} \right),
\end{align}
where
\begin{align}
\dot{\mathrm{U}}_{n+\alpha_m} &:= (1-\alpha_m) \dot{\mathrm{U}}_{n} + \alpha_m \dot{\mathrm{U}}_{n+1}, \\
\mathrm{U}_{n+\alpha_f} &:= (1-\alpha_f) \mathrm{U}_{n} + \alpha_f \mathrm{U}_{n+1},
\end{align}
$\Delta t_n := t_{n+1} - t_n$ is the time-step size, $t_{n+\alpha_f} := t_n + \alpha_f \Delta t_n$, and $\gamma$, $\alpha_m$, and $\alpha_f$ are free parameters.  The term $\dot{\mathrm{U}}_{n+\alpha_m}$ is often interpreted as an approximation of $\dot{\mathrm{U}}$ at time $t_{n+\alpha_m} := t_n + \alpha_m \Delta t_n$, while the term $\mathrm{U}_{n+\alpha_f}$ is often interpreted as an approximation of $\mathrm{U}$ at time $t_{n+\alpha_f}$.  The generalized-$\alpha$ method is second-order accurate if and only if
\begin{equation}
\gamma = \frac{1}{2} + \alpha_m - \alpha_f, \label{eq:gamma}
\end{equation}
and it is unconditionally stable if and only if
\begin{equation}
\alpha_m \geq \alpha_f \geq \frac{1}{2}.
\end{equation}
If Equation \eqref{eq:gamma} holds, then
\begin{align}
\mathrm{U}_{n+1} &= \mathrm{U}_{n} + \Delta t_n \left( \left(\frac{1}{2} - \alpha_m + \alpha_f\right) \dot{\mathrm{U}}_{n} + \left(\frac{1}{2} + \alpha_m - \alpha_f\right) \dot{\mathrm{U}}_{n+1} \right) \nonumber \\
&= \mathrm{U}_{n} + \Delta t_n \left( (1-\alpha_m) \dot{\mathrm{U}}_{n} + \alpha_m \dot{\mathrm{U}}_{n+1} + \left(\alpha_f - \frac{1}{2}\right) \dot{\mathrm{U}}_{n} - \left(\alpha_f - \frac{1}{2}\right) \dot{\mathrm{U}}_{n+1} \right) \nonumber \\
&= \mathrm{U}_{n} + \Delta t_n \left( \dot{\mathrm{U}}_{n+\alpha_m}+ \left(\alpha_f - \frac{1}{2}\right) \dot{\mathrm{U}}_{n} - \left(\alpha_f - \frac{1}{2}\right) \dot{\mathrm{U}}_{n+1} \right).
\end{align}
Thus, if the generalized-$\alpha$ method is second order-accurate, then
\begin{align}
\dot{\mathrm{U}}_{n+\alpha_m} &= \frac{\mathrm{U}^+_{n+\alpha_f} - \mathrm{U}^-_{n+\alpha_f}}{\Delta t_n} \label{eq:central_diff}
\end{align}
where
\begin{align}
\mathrm{U}^+_{n+\alpha_f} &:= \mathrm{U}_{n+1} + \left(\alpha_f - \frac{1}{2}\right) \Delta t_n \dot{\mathrm{U}}_{n+1}
\end{align}
and
\begin{align}
\mathrm{U}^-_{n+\alpha_f} &:= \mathrm{U}_{n} + \left(\alpha_f - \frac{1}{2}\right) \Delta t_n \dot{\mathrm{U}}_{n}.
\end{align}
Note that $\mathrm{U}^+_{n+\alpha_f}$ may be viewed as an approximation of $\mathrm{U}$ at time $t^+_{n+\alpha_f} := t_{n+1} + \left(\alpha_f-1/2\right) \Delta t_n$ due to the Taylor series
\begin{align}
\mathrm{U}\left(t^+_{n+\alpha_f}\right) &= \mathrm{U}(t_{n+1}) + \left(\alpha_f - \frac{1}{2}\right) \Delta t_n \dot{\mathrm{U}}(t_{n+1}) + O(\Delta t_n^2),
\end{align}
while $\mathrm{U}^-_{n+\alpha_f}$ may be viewed as an approximation of $\mathrm{U}$ at time $t^-_{n+\alpha_f} := t_n + \left(\alpha_f-1/2\right) \Delta t_n$ due to the Taylor series
\begin{align}
\mathrm{U}\left(t^-_{n+\alpha_f} \right) &= \mathrm{U}(t_{n}) + \left(\alpha_f - \frac{1}{2}\right) \Delta t_n \dot{\mathrm{U}}(t_{n}) + O(\Delta t_n^2).
\end{align}
Consequently, while $\dot{\mathrm{U}}_{n+\alpha_m}$ is often interpreted as an approximation of $\dot{\mathrm{U}}$ at time $t_{n+\alpha_m}$, Equation \eqref{eq:central_diff} indicates it may be instead be seen as a central difference approximation of $\dot{\mathrm{U}}$ at time $t_{n+\alpha_f}$ when the generalized-$\alpha$ method is second-order accurate.

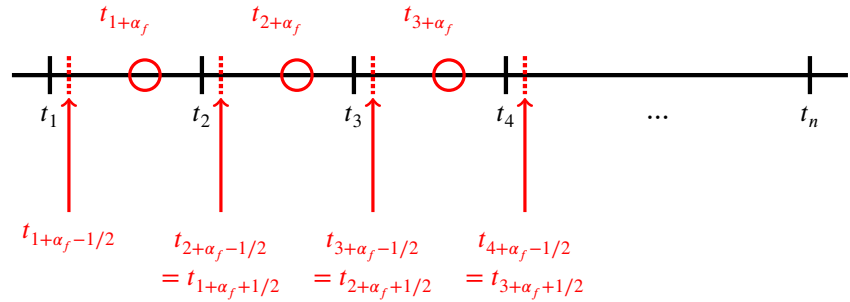
\begin{figure}[t!]
\begin{center}
% setting height of nodes
\newlength{\smallestHeight}
\setlength{\smallestHeight}{11pt}

\tikzset{block2/.style={shape=rectangle, draw, node distance=-1pt, minimum width = \stackWidth, line width=1pt, inner sep=0pt}}

\begin{tikzpicture}

% main line
\draw[ultra thick] (-0.5,0)--(10.5,0);

% tick marks and labels
\foreach \x/\y in {0/t_1, 2/t_2, 4/t_3, 6/t_4,  10/t_n}
 \draw[ultra thick] (\x,0.25) -- (\x,-0.25) node[below]{$\y$};
 % dots
\draw (8,0) node[below=4mm]{$\cdots$};

% intermediate points and labels
\foreach \i in {1.25,3.25, 5.25}% points on line
      \fill[very thick,draw=red,fill=none] (\i,0) circle (2 mm);
\draw[red] (1,0) node[above=4mm]{$t_{1+\alpha_f}$};
\draw[red] (3,0) node[above=4mm]{$t_{2+\alpha_f}$};
\draw[red] (5,0) node[above=4mm]{$t_{3+\alpha_f}$};

% TS points
\foreach \x in {0.25, 2.25, 4.25, 6.25}
 \draw[red, ultra thick, densely dotted] (\x,0.25) -- (\x,-0.25);% node[below]{$\y$};
 
 % TS labels
 \draw[<-, very thick, red] (0.25,-0.3)--++(-90:1.5) node[minimum height=\smallestHeight,below]{$t_{1+\alpha_f -1/2}$};
 \draw[<-, very thick, red] (2.25,-0.3)--++(-90:1.5) node[minimum height=\smallestHeight,below](first){$\begin{array}{c} t_{2+\alpha_f -1/2} \\ =t_{1+\alpha_f +1/2} \end{array}$};
 \draw[<-, very thick, red] (4.25,-0.3)--++(-90:1.5) node[minimum height=\smallestHeight,below](first){$\begin{array}{c} t_{3+\alpha_f -1/2} \\ =t_{2+\alpha_f +1/2} \end{array}$};
 \draw[<-, very thick, red] (6.25,-0.3)--++(-90:1.5) node[minimum height=\smallestHeight,below](first){$\begin{array}{c} t_{4+\alpha_f -1/2} \\ =t_{3+\alpha_f +1/2} \end{array}$};
 
 %{$\substack{t_{1+\alpha_f +1/2}\\= t_{2+\alpha_f -1/2}}$};

\end{tikzpicture}
\end{center}
\caption{Visual representation of the shifted temporal mesh on which the generalized-$\alpha$ method may be interpreted as an implicit midpoint method when it is second-order accurate and the original temporal mesh is uniform.}
\label{fig:timeMesh}
\end{figure}

For a non-uniform temporal mesh, we generally have that $U^+_{n+\alpha_f} \neq U^-_{(n+1)+\alpha_f}$ and $t^+_{n+\alpha_f} \neq t^-_{(n+1)+\alpha_f}$.  However, for a uniform temporal mesh, we have $U^+_{n+\alpha_f} = U^-_{(n+1)+\alpha_f}$ and $t^+_{n+\alpha_f} = t^-_{(n+1)+\alpha_f}$ for all $n \in \mathbb{N}$.  Defining in this case
\begin{align}
\mathrm{U}_{n+\alpha_f-1/2} &:= \mathrm{U}_{n} + \left(\alpha_f - \frac{1}{2}\right) \Delta t \dot{\mathrm{U}}_{n}
\end{align}
for $n \in \mathbb{N}$ where $\Delta t$ is the uniform time-step size, we have that
\begin{align}
\dot{\mathrm{U}}_{n+\alpha_m} &= \frac{\mathrm{U}_{n+\alpha_f+1/2} - \mathrm{U}_{n+\alpha_f-1/2}}{\Delta t} \label{eq:gold}
\end{align}
for $n \in \mathbb{N}$ when the generalized-$\alpha$ method is second order-accurate.  The terms $\left\{ \mathrm{U}_{n+\alpha_f-1/2} \right\}_{n \in \mathbb{N}}$ may be viewed as approximations of $\mathrm{U}$ on the shifted temporal mesh $\left\{ t_{n+\alpha_f-\frac{1}{2}} \right\}_{n \in \mathbb{N}}$ where
\begin{align}
t_{n+\alpha_f-\frac{1}{2}} := t_n + \left(\alpha_f - \frac{1}{2}\right) \Delta t
\end{align}
for $n \in \mathbb{N}$, and it can be shown that  $\left\{ \mathrm{U}_{n+\alpha_f-1/2} \right\}_{n \in \mathbb{N}}$ are in fact second-order approximations of $\left\{ \mathrm{U}\left(t_{n+\alpha_f-1/2}\right) \right\}_{n \in \mathbb{N}}$ when $\left\{ \mathrm{U}_{n} \right\}_{n \in \mathbb{N}}$ are second-order approximations of $\left\{ \mathrm{U}\left(t_{n}\right) \right\}_{n \in \mathbb{N}}$ and $\left\{ \dot{\mathrm{U}}_{n} \right\}_{n \in \mathbb{N}}$ are first-order approximations of $\left\{ \dot{\mathrm{U}}\left(t_{n}\right) \right\}_{n \in \mathbb{N}}$. Thus, when the generalized-$\alpha$ method is second-order accurate and the temporal mesh is uniform, it may be viewed as an implicit midpoint method on a shifted temporal mesh (see Figure \ref{fig:timeMesh}).  One might expect, then, the generalized-$\alpha$ method to inherit the conservation properties of the implicit midpoint method.  We demonstrate later this is indeed the case.

%fffffffffffffffff---------------------------------------------------------------------------------
%fffffffffffffffff--

In practice, the generalized-$\alpha$ parameters are typically chosen to be equal to
\begin{align}
\alpha_m &= \frac{1}{2} \left( \frac{3-\rho_{\infty}}{1+\rho_{\infty}} \right), \label{eq:alpha_m} \\
\alpha_f &= \frac{1}{1+\rho_{\infty}}, \label{eq:alpha_f}
\end{align}
where $\rho_{\infty} \in [0,1]$.  For this choice of parameters, the generalized-$\alpha$ method exhibits an optimal combination of accuracy in the low-frequency range and numerical damping in the high-frequency range for a linear model problem\cite{jansen2000generalized}.  The parameter $\rho_{\infty}$ then corresponds to the spectral radius of the amplification at infinite time step, and a choice of $\rho_{\infty} = 0$ annihilates the highest frequency in one step while a choice of $\rho_{\infty} = 1$ preserves the highest frequency.  As $\rho_{\infty}$ shifts from $\rho_{\infty} = 1$ to $\rho_{\infty} = 0$, $\alpha_f$ shifts from $\alpha_f = \frac{1}{2}$ to $\alpha_f = 1$, and as $\alpha_f$ shifts closer to $\alpha_f = 1$, the generalized-$\alpha$ method exhibits increased numerical damping.  In light of our interpretation of the generalized-$\alpha$ method as an implicit midpoint method on a shifted temporal mesh, we thus can view the generalized-$\alpha$ method with $\alpha_f > \frac{1}{2}$ as ``upwinding in time''.

If Equations \eqref{eq:alpha_m} and \eqref{eq:alpha_f} hold, it can be shown that second-order accuracy dictates $\gamma = \alpha_f$.  It follows then that the governing equations of the generalized-$\alpha$ method can be simply written entirely in terms of $\alpha_f$.  In particular, if the generalized-$\alpha$ method is second-order accurate, the temporal mesh is uniform, and Equations \eqref{eq:alpha_m} and \eqref{eq:alpha_f} hold, then the governing equations at the $n^\text{th}$ time step are
\begin{align}
&\mathrm{R}\left(\frac{\mathrm{U}_{n+\alpha_f+1/2} - \mathrm{U}_{n+\alpha_f-1/2}}{\Delta t},\mathrm{U}_{n+\alpha_f},t_{n+\alpha_f}\right) = \mathrm{0}, \\
&\mathrm{U}_{n+1} = \mathrm{U}_{n} + \Delta t \left( (1-\alpha_f) \dot{\mathrm{U}}_{n} + \alpha_f \dot{\mathrm{U}}_{n+1} \right), \\
&\mathrm{U}_{n+\alpha_f} = (1-\alpha_f) \mathrm{U}_{n} + \alpha_f \mathrm{U}_{n+1}, \\
&\mathrm{U}_{n+\alpha_f+1/2} = \mathrm{U}_{n+1} + \left(\alpha_f - \frac{1}{2}\right) \Delta t \dot{\mathrm{U}}_{n+1}, \\
&\mathrm{U}_{n+\alpha_f-1/2} = \mathrm{U}_{n} + \left(\alpha_f - \frac{1}{2}\right) \Delta t \dot{\mathrm{U}}_{n},
\end{align}
and $\alpha_f \in \left[\frac{1}{2},1\right]$ with $\alpha_f = \frac{1}{2}$ corresponding to no damping (and the implicit midpoint method on the original temporal mesh) and $\alpha_f = 1$ corresponding to maximal damping (and an implicit midpoint method on a temporal mesh shifted to the right by $\frac{\Delta t}{2}$).

\section{Application to the advection-diffusion problem}\label{sec:advection_diffusion}

Now consider the following advection-diffusion problem: Find $u: \Omega \times [0,\infty) \rightarrow \mathbb{R}$ such that
\begin{align}
\frac{\partial u}{\partial t} + \nabla \cdot \left( \mathbf{a} u - \kappa \nabla u \right) &= f && \text{ in } \Omega \times (0,\infty) \label{eq:ad_strong} \\
\kappa \nabla u \cdot \mathbf{n} - \min(\mathbf{a} \cdot \mathbf{n},0) u &= h && \text{ on } \Gamma \times (0,\infty) \\
u(\cdot,0) &= u_0 && \text{ in } \Omega
\end{align}
where $\Omega \subset \mathbb{R}^d$ is a $d$-dimensional spatial domain with $d \in \mathbb{N}$, $\Gamma$ is the boundary of $\Omega$, $\mathbf{n}$ is the outward unit normal vector to $\Omega$, $\mathbf{a}: \Omega \times (0,\infty) \rightarrow \mathbb{R}^d$ is the advection velocity satisfying $\nabla \cdot \mathbf{a} \equiv 0$, $\kappa : \Omega \times (0,\infty) \rightarrow \mathbb{R}$ is the diffusivity satisfying $\kappa > 0$, $f : \Omega \times (0,\infty) \rightarrow \mathbb{R}$ is the applied body force, $h : \Gamma \times (0,\infty) \rightarrow \mathbb{R}$ is the applied flux, and $u_0: \Omega \rightarrow \mathbb{R}$ is the applied initial condition.  Note that over the inflow boundary $\Gamma_\text{in} := \left\{ \mathbf{x} \in \Gamma: \mathbf{a}(\mathbf{x}) \cdot \mathbf{n} (\mathbf{x}) < 0 \right\}$, the sum of diffusive and advective fluxes is specified, while along the outflow boundary $\Gamma_\text{out} := \left\{ \mathbf{x} \in \Gamma: \mathbf{a}(\mathbf{x}) \cdot \mathbf{n} (\mathbf{x}) \geq 0 \right\}$, only the diffusive flux is specified.  This is necessary to arrive at a well-posed problem \cite{moghadam2011comparison}.  By integrating Equation \eqref{eq:ad_strong} over the spatial domain and invoking the divergence theorem, we attain
\begin{align}
\frac{d}{dt} \int_{\Omega} u d\Omega &= \int_{\Omega} f d\Omega + \int_{\Gamma} h d\Gamma - \int_{\Gamma_\text{out}} \left( \mathbf{a} \cdot \mathbf{n} \right) u d\Gamma
\end{align}
for all $t \in (0,\infty)$, and by integrating between times $t_\text{begin}$ and $t_\text{end}$ with $0 \leq t_\text{begin} \leq t_\text{end}$ and invoking the fundamental theorem of calculus, we further attain
\begin{align}
\int_{\Omega} u(\cdot,t_\text{end}) d\Omega &= \int_{\Omega} u(\cdot,t_\text{begin}) d\Omega + \int_{t_\text{begin}}^{t_\text{end}} \left( \int_{\Omega} f d\Omega + \int_{\Gamma} h d\Gamma - \int_{\Gamma_\text{out}} \left( \mathbf{a} \cdot \mathbf{n} \right) u d\Gamma \right) dt. \label{eq:ad_conservation}
\end{align}
The above is an integral balance law that we typically wish to preserve in some sense at the discrete level.

A stabilized or unstabilized Galerkin semi-discretization of the considered advection-diffusion problem takes the form: Find $u^h(t) \in \mathcal{V}^h$ for all $t \in [0,\infty)$ such that
\begin{align}
\int_{\Omega} \frac{\partial u^h}{\partial t} w^h d\Omega - \int_{\Omega} \left( \mathbf{a} u^h - \kappa \nabla u^h \right) \cdot \nabla w^h d\Omega + \int_{\Gamma_\text{out}} \left( \mathbf{a} \cdot \mathbf{n} \right) u^h w^h d\Gamma + S^h(u^h,w^h) &= \int_{\Omega} f w^h d\Omega + \int_{\Gamma} h w^h d\Gamma \label{eq:ad_galerkin}
\end{align}
for all $w^h \in \mathcal{V}^h$ and $t \in (0,\infty)$ and
\begin{align}
\int_{\Omega} u^h(\cdot,0) w^h d\Omega &= \int_{\Omega} u_0 w^h d\Omega \label{eq:ad_galerkin_ic}
\end{align}
for all $w^h \in \mathcal{V}^h$ where $\mathcal{V}^h$ is a finite-dimensional subspace of $H^1(\Omega)$ and $S^h : \mathcal{V}^h \times \mathcal{V}^h \rightarrow \mathbb{R}$ is a stabilization form ($S^h \equiv 0$ when no stabilization is applied).  Provided that $1 \in \mathcal{V}^h$ and $S^h(v^h,1) = 0$ for all $v^h \in \mathcal{V}^h$, we can take $w^h \equiv 1$ in Equation \eqref{eq:ad_galerkin}, integrate between times $t_\text{begin}$ and $t_\text{end}$ with $0 \leq t_\text{begin} \leq t_\text{end}$, and invoke the fundamental theorem of calculus to arrive at
\begin{align}
\int_{\Omega} u^h(\cdot,t_\text{end}) d\Omega &= \int_{\Omega} u^h(\cdot,t_\text{begin}) d\Omega + \int_{t_\text{begin}}^{t_\text{end}} \left( \int_{\Omega} f d\Omega + \int_{\Gamma} h d\Gamma - \int_{\Gamma_\text{out}} \left( \mathbf{a} \cdot \mathbf{n} \right) u^h d\Gamma \right) dt. \label{eq:ad_conservation_semi}
\end{align}
Thus a Galerkin semi-discretization inherits the integral balance law given in Equation \eqref{eq:ad_conservation} provided that $1 \in \mathcal{V}^h$ and $S^h(v^h,1) = 0$ for all $v^h \in \mathcal{V}^h$.  The property that $1 \in \mathcal{V}^h$ holds for most finite element approximation spaces that are used in practice.  The property that $S^h(v^h,1) = 0$ for all $v^h \in \mathcal{V}^h$ also holds for most stabilization methodologies that are used in practice.  In particular, it holds for the popular Streamline Upwind Petrov Galerkin (SUPG) method \cite{hughes1987recent}, the Galerkin Least Squares (GLS) method \cite{shakib1991new}, the Variational Multiscale (VMS) method \cite{bazilevs2007variational}, the method of orthogonal subscales \cite{codina2002stabilized}, edge stabilization \cite{burman2007continuous}, and local projection stabilization \cite{braack2006local}.  It also holds when viscosity-based discontinuity capturing operators \cite{bazilevs2007yzbeta} are employed.

The generalized-$\alpha$ method can be employed to discretize the first-order initial-value problem given by Equations \eqref{eq:ad_galerkin} and \eqref{eq:ad_galerkin_ic}.  This gives rise to the following governing equations at the $n^\text{th}$ time step:
\begin{align}
&\int_{\Omega} \dot{u}^h_{n+\alpha_m} w^h d\Omega - \int_{\Omega} \left( \mathbf{a} u^h_{n+\alpha_f} - \kappa_{n+\alpha_f} \nabla u^h_{n+\alpha_f} \right) \cdot \nabla w^h d\Omega + \int_{\Gamma_\text{out}} \left( \mathbf{a}_{n+\alpha_f} \cdot \mathbf{n} \right) u^h_{n+\alpha_f} w^h d\Gamma + S^h(u^h_{n+\alpha_f},w^h) \\
&= \int_{\Omega} f_{n+\alpha_f} w^h d\Omega + \int_{\Gamma} h_{n+\alpha_f} w^h d\Gamma
\end{align}
and
\begin{align}
u^h_{n+1} = u^h_{n} + \Delta t_n \left( (1-\gamma) \dot{u}^h_{n} + \gamma \dot{u}^h_{n+1} \right)
\end{align}
where
\begin{align}
\dot{u}^h_{n+\alpha_m} &:= (1-\alpha_m) \dot{u}^h_{n} + \alpha_m \dot{u}^h_{n+1}, \\
u^h_{n+\alpha_f} &:= (1-\alpha_f) u^h_{n} + \alpha_f u^h_{n+1},
\end{align}
$u^h_n$ and $u^h_{n+1}$ are the approximations of $u^h$ at times $t_n$ and $t_{n+1}$, $\dot{u}^h_n$ and $\dot{u}^h_{n+1}$ are the approximations of $\frac{\partial u^h}{\partial t}$ at times $t_n$ and $t_{n+1}$, $\kappa_{n+\alpha_f} = \kappa(\cdot,t_{n+\alpha_f})$, $\mathbf{a}_{n+\alpha_f} = \mathbf{a}(\cdot,t_{n+\alpha_f})$, $f_{n+\alpha_f} = f(\cdot,t_{n+\alpha_f})$, and $h_{n+\alpha_f} = h(\cdot,t_{n+\alpha_f})$.  If the generalized-$\alpha$ method is second-order accurate and the temporal mesh is uniform, Equation \eqref{eq:gold} applies and we can write
\begin{align}
\dot{u}^h_{n+\alpha_m} &= \frac{u^h_{n+\alpha_f+1/2} - \mathrm{u^h}_{n+\alpha_f-1/2}}{\Delta t}
\end{align}
where
\begin{align}
u^h_{n+\alpha_f+1/2} &:= u^h_{n+1} + \left(\alpha_f - \frac{1}{2}\right) \Delta t \dot{u}^h_{n+1}, \\
u^h_{n+\alpha_f-1/2} &:= u^h_{n} + \left(\alpha_f - \frac{1}{2}\right) \Delta t \dot{u}^h_{n}.
\end{align}
In this case, if $1 \in \mathcal{V}^h$ and $S^h(v^h,1) = 0$ for all $v^h \in \mathcal{V}^h$, we can take $w^h \equiv 1$ to immediately arrive at
\begin{align}
\int_{\Omega} u^h_{n+\alpha_f+1/2} d\Omega = \int_{\Omega} u^h_{n+\alpha_f-1/2} d\Omega + \Delta t \left( \int_{\Omega} f_{n+\alpha_f} d\Omega + \int_{\Gamma} h_{n+\alpha_f} d\Gamma - \int_{\Gamma_\text{out}} \left( \mathbf{a}_{n+\alpha_f} \cdot \mathbf{n} \right) u^h_{n+\alpha_f} d\Gamma \right).
\end{align}
We can sum over time steps to arrive at the balance law
\begin{align}
\int_{\Omega} u^h_{n_{\text{end}}+\alpha_f-1/2} d\Omega = \int_{\Omega} u^h_{n_{\text{begin}}+\alpha_f-1/2} d\Omega + \sum_{n=n_{\text{begin}}}^{n_{\text{end}}-1} \Delta t \left( \int_{\Omega} f_{n+\alpha_f} d\Omega + \int_{\Gamma} h_{n+\alpha_f} d\Gamma - \int_{\Gamma_\text{out}} \left( \mathbf{a}_{n+\alpha_f} \cdot \mathbf{n} \right) u^h_{n+\alpha_f} d\Gamma \right) \label{eq:ad_conservation_fully}
\end{align}
for two integers $1 \leq n_{\text{begin}} \leq n_{\text{end}}$.  The above is a fully-discrete analogue of Equation \eqref{eq:ad_conservation_semi}.  Thus, application of the generalized-$\alpha$ method to a conservative Galerkin semi-discretization of the advection-diffusion equation yields a conservative fully-discrete method if the generalized-$\alpha$ method is second-order accurate and the temporal mesh is uniform.  Similar results to those seen here can be attained if Dirichlet boundary conditions are applied provided a Lagrange multiplier field is introduced\cite{evans2013isogeometric}, and local conservation results can also be attained using the method described by Hughes et. al\cite{hughes2000continuous}.

\section{Application to systems of conservation laws}\label{sec:conservation_laws}

We last consider a general system of differential conservation laws.  Without loss of generality, we consider only periodic boundary conditions.  The problem of interest then reads as follows: Find $\textbf{U} : \Omega \times [0,\infty) \rightarrow \mathbb{R}^p$ such that
\begin{align}
\frac{\partial \mathbf{U}}{\partial t} + \sum_{i=1}^{d} \frac{\partial \mathbf{F}^i}{\partial x_i} &= \mathbf{S} && \text{ in } \Omega \times (0,\infty) \label{eq:cons} \\
\mathbf{U}(\cdot,0) &= \mathbf{U}_0 && \text{ in } \Omega
\end{align}
and $\textbf{U}$ is periodic in each spatial direction where $\Omega := (0,1)^d$ is the $d$-dimensional domain with $d \in \mathbb{N}$, $\left\{ \textbf{F}^i \right\}_{i=1}^d$ are the flux vectors in each spatial direction, $\textbf{S}$ is the source vector, and $\mathbf{U}_0 : \Omega \rightarrow \mathbb{R}^p$ is the applied initial condition.  The flux vector and source vectors can depend on space, time, as well as $\textbf{U}$ and its spatial derivatives in each direction.  We require, however, that the flux vector and source vectors are themselves periodic in each direction.  Both the Euler and Navier-Stokes equations can be written in the above form, as can the equations governing magnetohydrodynamics.  Integrating \eqref{eq:cons} over the spatial domain, invoking the divergence theorem, and then integrating between times $t_\text{begin}$ and $t_\text{end}$ with $0 \leq t_\text{begin} \leq t_\text{end}$ and invoking the fundamental theorem of calculus, we attain
\begin{align}
\int_{\Omega} \mathbf{U}(\cdot,t_\text{end}) d\Omega = \int_{\Omega} \mathbf{U}(\cdot,t_\text{begin}) d\Omega + \int_{t_\text{begin}}^{t_\text{end}} \int_{\Omega} \mathbf{S} d\Omega dt \label{eq:cons_conservation}
\end{align}
which is a generalization of Equation \eqref{eq:ad_conservation} to the current setting.  A stabilized or unstabilized Galerkin semi-discretization of the above problem using conservation variables takes the form: Find $\mathbf{U}^h(t) \in \bm{\mathcal{V}}^h$ for all $t \in [0,\infty)$ such that
\begin{align}
\int_{\Omega} \frac{\partial \mathbf{U}^h}{\partial t} \cdot \mathbf{W}^h d\Omega - \sum_{i=1}^d \int_{\Omega} \mathbf{F}^i \cdot \frac{\partial \mathbf{W}^h}{\partial x_i} d\Omega + S^h\left(\mathbf{U}^h,\mathbf{W}^h\right) &= \int_{\Omega} \mathbf{S} \cdot \mathbf{W}^h d\Omega \label{eq:cons_galerkin}
\end{align}
for all $\mathbf{W}^h \in \bm{\mathcal{V}}^h$ and $t \in (0,\infty)$ and
\begin{align}
\int_{\Omega} \mathbf{U}^h(\cdot,0) \cdot \mathbf{W}^h d\Omega &= \int_{\Omega} \mathbf{U}_0 \cdot \mathbf{W}^h  d\Omega \label{eq:cons_galerkin_ic}
\end{align}
for all $\mathbf{W}^h \in \bm{\mathcal{V}}^h$ where $\bm{\mathcal{V}}^h$ is a finite-dimensional subspace of $(H^1_\text{per}(\Omega))^p$ and $S^h : \bm{\mathcal{V}}^h \times \bm{\mathcal{V}}^h \rightarrow \mathbb{R}$ is a stabilization form (again $S^h \equiv 0$ in the absence of stabilization).  If the unit vector $\mathbf{e}_j$ is a member of $\bm{\mathcal{V}}^h$ and $S^h(\mathbf{v}^h,\mathbf{e}_j)$ for all $\mathbf{v}^h \in \bm{\mathcal{V}}^h$ for each $j = 1, \ldots, p$, then we can use the same procedure employed to arrive at Equation \eqref{eq:ad_conservation_semi} to also arrive at
\begin{align}
\int_{\Omega} \mathbf{U}^h(\cdot,t_\text{end}) d\Omega = \int_{\Omega} \mathbf{U}^h(\cdot,t_\text{begin}) d\Omega + \int_{t_\text{begin}}^{t_\text{end}} \int_{\Omega} \mathbf{S} d\Omega dt \label{eq:cons_conservation_semi}
\end{align}
for $0 \leq t_\text{begin} \leq t_\text{end}$.  Application of the generalized-$\alpha$ method to the Galerkin semi-discretization results in the following governing equations at the $n^\text{th}$ time step:
\begin{align}
\int_{\Omega} \dot{\mathbf{U}}_{n+\alpha_m}^h \cdot \mathbf{W}^h d\Omega - \sum_{i=1}^d \int_{\Omega} \mathbf{F}_{n+\alpha_f}^i \cdot \frac{\partial \mathbf{W}^h}{\partial x_i} d\Omega + S^h\left(\mathbf{U}_{n+\alpha_f}^h,\mathbf{W}^h\right) &= \int_{\Omega} \mathbf{S}_{n+\alpha_f} \cdot \mathbf{W}^h d\Omega
\end{align}
and
\begin{align}
\mathbf{U}^h_{n+1} = \mathbf{U}^h_{n} + \Delta t_n \left( (1-\gamma) \dot{\mathbf{U}}^h_{n} + \gamma \dot{\mathbf{U}}^h_{n+1} \right)
\end{align}
where
\begin{align}
\dot{\mathbf{U}}^h_{n+\alpha_m} &:= (1-\alpha_m) \dot{\mathbf{U}}^h_{n} + \alpha_m \dot{\mathbf{U}}^h_{n+1}, \\
\mathbf{U}^h_{n+\alpha_f} &:= (1-\alpha_f) \mathbf{U}^h_{n} + \alpha_f \mathbf{U}^h_{n+1},
\end{align}
$\mathbf{U}^h_n$ and $\mathbf{U}^h_{n+1}$ are the approximations of $\mathbf{U}^h$ at times $t_n$ and $t_{n+1}$, $\dot{\mathbf{U}}^h_n$ and $\dot{\mathbf{U}}^h_{n+1}$ are the approximations of $\frac{\partial \mathbf{U}^h}{\partial t}$ at times $t_n$ and $t_{n+1}$, $\left\{ \mathbf{F}^i_{n+\alpha_f} \right\}_{i=1}^d$ are the flux vectors $\left\{ \mathbf{F}^i \right\}_{i=1}^d$ evaluated at time $t_{n+\alpha_f}$ using the value and spatial derivatives of $\mathbf{U}^h_{n+\alpha_f}$, and $\mathbf{S}_{n+\alpha_f}$ is the source vector evaluated at time $t_{n+\alpha_f}$ using the value and spatial derivatives of $\mathbf{U}^h_{n+\alpha_f}$.  If the generalized-$\alpha$ method is second-order accurate and the temporal mesh is uniform, then
\begin{align}
\dot{\mathbf{U}}^h_{n+\alpha_m} &= \frac{\mathbf{U}^h_{n+\alpha_f+1/2} - \mathrm{\mathbf{U}^h}_{n+\alpha_f-1/2}}{\Delta t}
\end{align}
where
\begin{align}
\mathbf{U}^h_{n+\alpha_f+1/2} &:= \mathbf{U}^h_{n+1} + \left(\alpha_f - \frac{1}{2}\right) \Delta t \dot{\mathbf{U}}^h_{n+1}, \\
\mathbf{U}^h_{n+\alpha_f-1/2} &:= \mathbf{U}^h_{n} + \left(\alpha_f - \frac{1}{2}\right) \Delta t \dot{\mathbf{U}}^h_{n},
\end{align}
and if the unit vector $\mathbf{e}_j$ is a member of $\bm{\mathcal{V}}^h$ and $S^h(\mathbf{v}^h,\mathbf{e}_j)$ for all $\mathbf{v}^h \in \bm{\mathcal{V}}^h$ for each $j = 1, \ldots, p$, then we can use the same procedure employed to arrive at Equation \eqref{eq:ad_conservation_fully} to also arrive at the following fully-dicrete analogue of Equation \eqref{eq:cons_conservation_semi}:
\begin{align}
\int_{\Omega} \mathbf{U}^h_{n_{\text{end}}+\alpha_f-1/2} d\Omega = \int_{\Omega} \mathbf{U}^h_{n_{\text{begin}}+\alpha_f-1/2} d\Omega + \sum_{n=n_{\text{begin}}}^{n_{\text{end}}-1} \Delta t \int_{\Omega} \mathbf{S}_{n+\alpha_f} d\Omega \label{eq:cons_conservation_fully}
\end{align}
for two integers $1 \leq n_{\text{begin}} \leq n_{\text{end}}$.  Therefore we also have that application of the generalized-$\alpha$ method to a conservative Galerkin semi-discretization of a system of differential conservation laws yields a conservative fully-discrete method if the generalized-$\alpha$ method is second-order accurate, the temporal mesh is uniform, and the conservation variables themselves are discretized.

\section{Discretization with nonconservation variables}\label{sec:nonconservative variables}

It is common practice to discretize systems of differential conservation laws using variables other than the conservation variables.  For instance, the use of pressure primitive variables or entropy variables is common in the discretization of the Euler and Navier-Stokes equations\cite{hughes1986new}.  A stabilized or unstabilized Galerkin semi-discretization of the system of differential conservation laws analyzed in the previous section using a set of nonconservation variables takes the form: Find $\mathbf{V}^h(t) \in \bm{\mathcal{V}}^h$ for all $t \in [0,\infty)$ such that
\begin{align}
\int_{\Omega} \left(\frac{\partial \mathbf{U}}{\partial \mathbf{V}}\left(\mathbf{V}^h\right) \frac{\partial \mathbf{V}^h}{\partial t}\right) \cdot \mathbf{W}^h d\Omega - \sum_{i=1}^d \int_{\Omega} \mathbf{F}^i \cdot \frac{\partial \mathbf{W}^h}{\partial x_i} d\Omega + S^h\left(\mathbf{U}(\mathbf{V}^h),\mathbf{W}^h\right) &= \int_{\Omega} \mathbf{S} \cdot \mathbf{W}^h d\Omega \label{eq:alt}
\end{align}
for all $\mathbf{W}^h \in \bm{\mathcal{V}}^h$ and $t \in (0,\infty)$ and
\begin{align}
\int_{\Omega} \mathbf{U}(\mathbf{V}^h(\cdot,0)) \cdot \mathbf{W}^h d\Omega &= \int_{\Omega} \mathbf{U}_0 \cdot \mathbf{W}^h  d\Omega \label{eq:alt_ic}
\end{align}
for all $\mathbf{W}^h \in \bm{\mathcal{V}}^h$ where $\mathbf{V}^h$ is the approximate vector of nonconservation variables and $\mathbf{U}: \mathbb{R}^p \rightarrow \mathbb{R}^p$ is the mapping between nonconservation variables and conservation variables.  The above Galerkin semi-discretization harbors the same conservation properties as a Galerkin semi-discretization using conservation variables.  The same is not true, however, for the fully-discrete method attained after application of the generalized-$\alpha$ method.  To see this, note that time-discretization of Equations \eqref{eq:alt} and \eqref{eq:alt_ic} using the generalized-$\alpha$ method results in the following governing equations at the $n^\text{th}$ time step:
\begin{align}
\int_{\Omega} \left(\frac{\partial \mathbf{U}}{\partial \mathbf{V}}\left(\mathbf{V}^h_{n+\alpha_f}\right) \dot{\mathbf{V}}_{n+\alpha_m}^h\right) \cdot \mathbf{W}^h d\Omega - \sum_{i=1}^d \int_{\Omega} \mathbf{F}_{n+\alpha_f}^i \cdot \frac{\partial \mathbf{W}^h}{\partial x_i} d\Omega + S^h\left(\mathbf{U}\left(\mathbf{V}_{n+\alpha_f}^h\right),\mathbf{W}^h\right) &= \int_{\Omega} \mathbf{S}_{n+\alpha_f} \cdot \mathbf{W}^h d\Omega \label{eq:cons_galerkin_non}
\end{align}
and
\begin{align}
\mathbf{V}^h_{n+1} = \mathbf{V}^h_{n} + \Delta t_n \left( (1-\gamma) \dot{\mathbf{V}}^h_{n} + \gamma \dot{\mathbf{V}}^h_{n+1} \right)
\end{align}
where
\begin{align}
\dot{\mathbf{V}}^h_{n+\alpha_m} &:= (1-\alpha_m) \dot{\mathbf{V}}^h_{n} + \alpha_m \dot{\mathbf{V}}^h_{n+1}, \\
\mathbf{V}^h_{n+\alpha_f} &:= (1-\alpha_f) \mathbf{V}^h_{n} + \alpha_f \mathbf{V}^h_{n+1}.
\end{align}
Unfortunately, even when the generalized-$\alpha$ method is second-order accurate and the temporal mesh is uniform, we cannot express $\frac{\partial \mathbf{U}}{\partial \mathbf{V}}\left(\mathbf{V}^h_{n+\alpha_f}\right) \dot{\mathbf{V}}_{n+\alpha_m}^h$ in terms of a difference of conservation variable states.  As such, we cannot arrive at a discrete balance law analogous to Equation \eqref{eq:cons_conservation_fully}.  We can remedy this situation by replacing $\frac{\partial \mathbf{U}}{\partial \mathbf{V}}\left(\mathbf{V}^h_{n+\alpha_f}\right) \dot{\mathbf{V}}_{n+\alpha_m}^h$ in Equation \eqref{eq:cons_galerkin_non} with
\begin{align}
\frac{\hat{\mathbf{U}}^h_{n+\alpha_f+1/2} - \mathrm{\hat{\mathbf{U}}^h}_{n+\alpha_f-1/2}}{\Delta t}
\end{align}
when Equation \eqref{eq:gamma} holds and the temporal mesh is uniform where
\begin{align}
\hat{\mathbf{U}}^h_{n+\alpha_f+1/2} &:= \mathbf{U}\left(\mathbf{V}^h_{n+1}\right) + \left(\alpha_f - \frac{1}{2}\right) \Delta t \frac{\partial \mathbf{U}}{\partial \mathbf{V}}\left(\mathbf{V}^h_{n+1}\right) \dot{\mathbf{V}}^h_{n+1}, \\
\hat{\mathbf{U}}^h_{n+\alpha_f-1/2} &:= \mathbf{U}\left(\mathbf{V}^h_{n}\right) + \left(\alpha_f - \frac{1}{2}\right) \Delta t \frac{\partial \mathbf{U}}{\partial \mathbf{V}}\left(\mathbf{V}^h_n\right) \dot{\mathbf{V}}^h_{n}.
\end{align}
The resulting fully-discrete method then admits the discrete balance law
\begin{align}
\int_{\Omega} \hat{\mathbf{U}}^h_{n_{\text{end}}+\alpha_f-1/2} d\Omega = \int_{\Omega} \hat{\mathbf{U}}^h_{n_{\text{begin}}+\alpha_f-1/2} d\Omega + \sum_{n=n_{\text{begin}}}^{n_{\text{end}}-1} \Delta t \int_{\Omega} \mathbf{S}_{n+\alpha_f} d\Omega
\end{align}
for two integers $1 \leq n_{\text{begin}} \leq n_{\text{end}}$. 

\section{Conclusion}\label{sec:conclusion}
In this short communication, we showed that application of the second-order accurate generalized-$\alpha$ method to a stabilized or unstabilized Galerkin semi-discretization of a system of differential conservation laws results in a fully-discrete method that inherits the conservation properties of the underlying semi-discretization provided the temporal mesh is uniform.  To do so, we first conducted a critical examination of the second-order accurate generalized-$\alpha$ method for first-order initial value problems, and we found it may be viewed as an implicit midpoint method on a shifted temporal mesh when the temporal mesh is uniform.  We then employed this knowledge to show second-order accurate generalized-$\alpha$ time integration of a conservative Galerkin semi-discretization of the advection-diffusion equation using  a uniform temporal mesh yields a fully-discrete method admitting a discrete balance law, and we then illustrated the same is true for general systems of differential conservation laws provided the conservation variables are themselves discretized.  When nonconservation variables are instead discretized, the resulting fully-discrete method is not conservative, but we demonstrated how to modify the generalized-$\alpha$ method to arrive at a discrete balance law in this case.  All the theoretical results appearing in this note have been verified by numerical experiments, but these experiments are not discussed here for brevity.

The theoretical results appearing in this note hold under the restrictive assumption of a uniform temporal mesh.  While a uniform temporal mesh is most commonly employed in practice, significant efficiency gains are sometimes possible with adaptive time integration.  We do not believe that discrete balance laws can be derived under the less restrictive assumption of a nonuniform temporal mesh, but we believe it is possible that the generalized-$\alpha$ method can be modified to ensure conservation in this case, just as was done for the case of nonconservation variables in this note.

Finally, while we only showed that the second-order accurate generalized-$\alpha$ method for first-order initial value problems can be viewed as an implicit midpoint method on a shifted temporal mesh if the unshifted temporal mesh is uniform, the same is true for second-order initial value problems as well.  To see this, note that application of the generalized-$\alpha$ method to a second-order initial value problem results in a residual equation of the form
\begin{align}
\mathrm{R}\left(\ddot{\mathrm{U}}_{n+\alpha_m},\dot{\mathrm{U}}_{n+\alpha_f},\mathrm{U}_{n+\alpha_f},t_{n+\alpha_f}\right) = \mathrm{0}
\end{align}
at the $n^\text{th}$ time step, and second-order accuracy still dictates that $\gamma = \frac{1}{2} + \alpha_m - \alpha_f$\cite{chung1994family}.  Consequently, the same analysis conducted in this note can also be used to show that
\begin{align}
\ddot{\mathrm{U}}_{n+\alpha_m} &= \frac{\dot{\mathrm{U}}_{n+\alpha_f+1/2} - \dot{\mathrm{U}}_{n+\alpha_f-1/2}}{\Delta t}
\end{align}
on a uniform temporal mesh where $\dot{\mathrm{U}}_{n+\alpha_f+1/2} := \dot{\mathrm{U}}_{n+1} + \left(\alpha_f - \frac{1}{2}\right) \Delta t \ddot{\mathrm{U}}_{n+1}$ and $\dot{\mathrm{U}}_{n+\alpha_f-1/2} = \dot{\mathrm{U}}_{n} + \left(\alpha_f - \frac{1}{2}\right) \Delta t \ddot{\mathrm{U}}_{n}$, and as a result, application of the second-order generalized-$\alpha$ method to a Galerkin elastodynamics semi-discretization results in a fully-discrete method with a discrete balance law for momentum.  We leave further analysis of this, as well as extension of the theoretical results shown here to higher-order generalizations of the generalized-$\alpha$ method\cite{behnoudfar2021higher,behnoudfar2021higher_parabolic}, for future work.

%\backmatter

\section*{Acknowledgments}
Both authors were partially funded by the Army Research Office under Award Number W911NF20P0002.

\subsection*{Author contributions}

Both authors contributed to the conceptualization, writing, and editing of this note.

\subsection*{Financial disclosure}

None reported.

\subsection*{Conflict of interest}

The authors declare no potential conflict of interests.

\section*{Supporting information}

There is no supporting information for this article.

\nocite{*}% Show all bib entries - both cited and uncited; comment this line to view only cited bib entries;

\bibliography{wileyNJD-AMA.bib}%

\end{document}